\documentclass[a4paper,authoryear,11pt]{article}

\usepackage{empheq} % Math symbols, fonts, environments
\usepackage{amsfonts,amsthm, amscd, epsfig, amsmath, amssymb,enumerate}
\usepackage{color}
\usepackage[T1]{fontenc}  
\usepackage{epstopdf}

\def \P {{\mathbb P}}
\def \E {{\mathbb E}}

\def \ra {\rightarrow }

\def \s {y}

\def\e{\epsilon}
\def\s{\sigma}

\def\l{\lambda}

\def\a{\alpha}

\begin{document}

\thispagestyle{empty}
\def\thefootnote{\fnsymbol{footnote}}

\title{Climb on the Bandwagon: Consensus and periodicity in a lifetime utility model with strategic interactions}

\author{Paolo Dai Pra\footnote{Department of Mathematics, University of Padova, 63, Via Trieste,  I - 35121 Padova, Italy;  daipra@math.unipd.it}, $\,$
Elena Sartori \footnote{Department of Mathematics, University of Padova, 63, Via Trieste,  I - 35121 Padova, Italy;  esartori@math.unipd.it}
$\,$ and $\,$ Marco Tolotti \footnote{Department of Management, Ca' Foscari University of Venice,
\mbox{S.Giobbe - Cannaregio 873, I - 30121 Venice, Italy;} tolotti@unive.it}
}

\maketitle{}

\begin{abstract}
What is the emergent long-run equilibrium of a society where many interacting agents bet on the optimal energy to put in place in order to \emph{climb on the Bandwagon}? In this paper we  study the collective behavior of a large population of agents being either Left or Right: the core idea is that agents benefit from being with the winner party, but, on the other hand, they suffer a cost in changing their status quo. {At the microscopic level the model is formulated as a stochastic, symmetric dynamic game with $N$ players. In the macroscopic limit as $N \rightarrow +\infty$, we obtain a {\em mean field game} whose equilibria describe the ``rational'' collective behavior of the society. It is of particular interest to detect  the emerging long-time attractors, e.g. consensus or oscillating behavior. Significantly, we discover that \emph{bandwagoning} can be persistent at the macro level: endogenously generated periodicity is in fact detected.}
%
%
%
%It is evident that there are two opposite drivers: climbing on the Bandwagon is desirable but costly. Transitions from Left to Right and vice versa are neither deterministic, nor a priori determined: agents themselves are in charge to choose the energy they want to put in place. Some questions naturally arise: (i) what is the \emph{optimal} level of energy agents are ready to spend in order to be with the majority? (ii) what are the emerging long-time attractors of the society? (consensus, oscillating behavior); %(iii) how to formalize the strategic interactions and the optimal control problem? 
%(iii) how to characterize the optimal trajectories (i.e., the Nash Equilibria) of this strategic game? (are they uniquely defined or not, smooth or oscillating). Significantly, we discover that \emph{bandwagoning} can be persistent at the macro level: endogenously generated periodicity is in fact detected.
\end{abstract}

\begin{quote}
\textbf{Keywords:}  Consensus, Mean field games; Multiple Nash Equilibria, Opinion dynamics; Social interactions
\end{quote}

%\maketitle

%%%%%%%%%%%%%%%%%%%%%%%%%%%%%%%%%%%%%%%%%%%%%%%%%%%%%%%%%%%%%%%%%%%%%%%%%%%%%%%%
\section*{Introduction}

The emergence of collective behavior in complex societies %, thought as the attitude of a society to express social norms to which single agents are willing (or forced) to conform, 
has been one of the most studied paradigms of social sciences in the last decades. Pioneering works, among the others, have been devoted to the study of segregation (see \cite{sc71}), social innovation (see \cite{sc73_helmets}), riot's formation (see \cite{Gr}), social distance (see \cite{ak97social}) or the emergence of  prices in financial markets  (see \cite{fo74}).

Following the celebrated \emph{micromotives and macrobehavior} paradigm by Shelling (see \cite{sh78_micro}), large attention has been paid to the mechanisms under which social norms emerge as an aggregate output of a large population of interacting agents.  In \cite{BeD}, we find one of the first attempts to formally analyze the large limits for economies characterized by social externalities, whereas in \cite{kalai04} a similar results is formalized in a game-theoretical setting, and in \cite{ep02} in the field of artificial societies and agent based models. One of the key factors underpinning all these models is the presence of \emph{positive externalities}, meaning that the single agent benefits from aligning with social norms, or, put differently, from being with the majority. The behavioral attitude behind this assumption is, basically, conformism, imitation or peer pressure. Of course, depending on the applications, different behavioral assumptions can be made (see, for instance, \cite{pnas16} or \cite{ch97} for nonconforming individuals or minority games, respectively). 

One aspect that, to our opinion, has not been sufficiently considered in the aforementioned literature, is the fact that changing opinion %to reach a majority or to meet the prevailing social norm, 
may be \emph{effort-demanding}. Gather information, modify habits or practices, revise operational strategies, join a new technology or, in one word, \emph{climbing on the bandwagon}, may result to be a costly operation.  This is the main goal of this paper: %build a tractable model where agents are willing to conform with the majority besides suffering a cost to change their mind. The main aim of this work is to 
provide a stylized model to describe a large population of conformist agents, in charge to optimally determine the level of energy required to stick with the majority.  Our idea is simple: the higher the effort is, the more likely it is to be with the winner party and the higher is the associated cost.

Indeed, to dynamically study the aggregate behavior of the society, we rely on a lifetime utility maximization problem where the agent is in charge to optimally set the effort to put in place to change status. Two remarks are needed. Firstly, the lifetime setup immediately refers to a parallel strand of literature, related to more classical optimization problems for consensus formation (see \cite{hu03}, \cite{Ieee13}, \cite{JOTA16}). The second remark pertains to the modeling structure of the society; in order to study into details (and possibly to obtain closed-form solutions) the relationship between social interaction and frictions in changing opinion, we stick with the simplest possible geometry: a mean-field model.\footnote{For recent literature investigating the relationship between the network geometry and the diffusion of knowledge, innovation, consensus, see \cite{ho06}, \cite{mo10}, \cite{yo11}, \cite{va12}.} This means, in particular, that aggregate statistics such as the empirical distribution (or the empirical mean) are sufficient to fully describe the Markovian system. With this respect, we move in the framework of {\em mean field games}.  The recent theory in this field has put forward a class of dynamic games for which the limit behavior, as the number of agents increases to infinity, can be described in analytic terms (see \cite{caines2005}, \cite{Li06}, \cite{Gu11}). % In this limit, the solution of the dynamic game is given by a system of two coupled equations: one is the Hamilton-Jacobi-Bellman equation for the value function, the second is the master equation for the optimal evolution of the representative agent. 
%Bandwagon (impact of digital bandwagon in elections \cite{PNAS15})  
 
In a nutshell, in our model the state of each agent (Left or Right) evolves as a controlled Markov process; the discounted lifetime utility is formed by two components: a reward received only when agreeing with the majority and a quadratic instantaneous cost related to the effort put in place. Since each agent looks for the control maximizing her own utility, under strategic interactions, it is natural to consider Nash equilibria for the resulting dynamic game. 
When the limit of infinitely many players is considered we formally obtain the usual {\em mean field game equation}, given by a system of two coupled equations: one is the Hamilton-Jacobi-Bellman equation for the value function, the second is the master equation for the optimal evolution of the representative agent. The rigorous foundation of this formal limit is to a large extent an open problem. Rigorous convergence results have been obtained recently for diffusion models  (see \cite{Go14},\cite{Ca15}) and for models with discrete state space (see \cite{CePe17}); this results, that are limited to the case of finite time horizon, require assumptions that guarantee the {\em uniqueness} of the solution of the mean field game equation. In the models considered here this uniqueness fails. As pointed out in \cite{CeFi17}, {\em all} solutions of the mean field game equation have ``physical'' significance for the $N$-player game: if the feedback control corresponding to {\em one solution} of the mean field game equation is applied by each player in the $N$-player game, it is an {\em approximate} Nash equilibrium, with the approximation error going to zero as $N \ra +\infty$. However, some solutions of the mean field game may not be obtained as limit of Nash equilibria of the $N$-player game. For examples of non uniqueness and the related convergence problems we refer to \cite{bafi17} and \cite{cedafipe18}.

In this paper we do not provide a rigorous convergence result in the number of players. We, rather, devote our attention to the long-run behavior of the asymptotic model and its properties.
 In particular, for the class of models introduced below, we find cases in which the long-run behavior of the mean field game leads to consensus, 
other in which the limit system admits periodic and non-constant  solutions. This rhythmic behavior, a sort of \emph{macro bandwagoning} of the society, emerges in absence of external periodic signals, and it is endogenously produced by the \emph{micro motives} behind the strategic behavior of agents. Investigation of periodic behavior of multi-agent systems has been recently studied, for instance, in \cite{yang2016} and \cite{basar2016}. In the context of mean-field games, however, periodic behavior
has been often predicted but, to our knowledge, proved only for the rather celebrated {\em Mexican wave} model (see \cite{Gu11}). It must be remarked that the Mexican wave model possesses a continuous symmetry, which allows the appearance of traveling wave solution. The model we propose below has a discrete (actually binary) space structure, so there is no continuous symmetry. Recent years have seen a formidable effort in the attempt of explaining rigorously the emergence of collective periodicity in noisy systems of interacting units. Given the impossibility of accounting for the huge related literature, we only mention the inspiring work \cite{Li04}, and few available rigorous results in \cite{tou14}, \cite{dp13}, \cite{dp15}, \cite{dit17}. In these works a key role in the emergence of periodicity is  played by {\em delay} in the information transmission (see \cite{tou14}, \cite{dit17}) and {\em dissipation}  (see  \cite{dp13}, \cite{dp15}). One of the main purposes of this paper is to show that collective periodic behavior can alsoresult from agents' utility optimization.

\section*{The microscopic model}

Consider a network of $N$ interacting agents, each possessing a binary state $\s_i(t) \in \{-1,1\}$ at time $t \geq 0$. 
Every agent can control her own state by means of the control $u_i = (u_i(t))_{t \geq 0}$. We assume here {\em close-loop} controls under complete information: 
\[
u_i(t) = \varphi_i(t, \s(t))
\]
for some function $\varphi_i$ which is right-continuous in $t$ and depending on the whole state $\s(t) = (\s_j(t))_{j=1}^N$ at time $t$.
The controlled stochastic dynamics  are given by
\begin{equation} \label{dyn}
\P\left(\s_i(t+h) = -\s_i(t) \, \Big| \,\s(s), \, s \leq t \right) = u_i(t) h + o(h).
\end{equation}
In other words, $u_i(t)$ is the probability rate of {\em flipping} the state $\s_i$.
Let 
\[
m_N(t) := \frac{1}{N} \sum_{i=1}^N \s_i(t)
\]
be the average state of the network at time $t$. 
%We introduce the {\em delay function} $\varphi(u) := \frac{k^{n+1}}{n!} u^n e^{-ku}$, which depends on the parameters $n \in  \Z^+$ and $k \in [0,+\infty]$. For $k = +\infty$, $\varphi$ should be meant as the Dirac distribution at $u=0$. 
The instantaneous reward of agent $i$ at time $t$ is given by
\[
R_i(t) :=  \s_i(t) m_N(t) - \frac{1}{2\mu(\s_i(t), m_N(t))} u^2_i(t).
\]
The two summands in the reward $R_i$ are easy to interpret. The term $\s_i(t) m_N(t)$ favors imitation: agents are conformist, they gain when aligned with  the majority. The term  $- \frac{\mu}{2} u^2_i(t)$ is an {\em energy} cost: a rapid change of the state would require high values for $u_i$, which are costly. The factor $\mu(\s_i(t), m_N(t))$, that we assume to be nonnegative,  modulates the relevance of this cost term: large values of $\mu$ allow high mobility to the agents, who can rapidly adapt to a change in the majority. Conversely, small values of $\mu$ reduce the adaptive response of agents. We allow $\mu$ to depend on the state of agent $i$ and on the average state of the network.
\\ Each agent $i$ aims at maximizing the discounted lifetime utility
\[
U_i := \E\left[\int_0^{+\infty} e^{-\l t} R_i(t) dt \right],
\]
where $\l>0$ is a constant discount factor.

A control $u^* = (u^*_1, u^*_2, \ldots, u^*_N)$ is called a {\em Nash equilibrium} if for every $i=1,\ldots,N$, assuming that all agents $j \neq i$ use the control $u_j^*$, we have $U_i(u^*_i) \geq U_i(u_i)$ for every other control $u_i$: in equilibrium no agent has interest in changing her strategy. Note that this dynamic game is invariant for permutation of agents, so it falls within the domain of {\em mean-field games}  (see \cite{Ca06}, \cite{Li06}).

\section*{The macroscopic model}

The limit as $N \rightarrow +\infty$ of the dynamic game described above is easy to obtain at a heuristic level. One expects that the average state $m_N(t)$ obeys a {\em Law of Large Numbers}, so it converges to a deterministic limit $m(t)$. The {\em representative agent} aims at maximizing
\begin{equation} \label{limmax}
J(u) := \E\left[\int_0^{+\infty} e^{-\l t} \left( \s(t) m(t) - \frac{1}{2\mu(\s(t), m(t))} u^2(t) \right) dt \right].
\end{equation}
An equilibrium control $u^*$ must satisfy the following consistency relation: if we denote by $\s^*(t)$ the process produced by the control $u^*$, then
\[
m(t) = \E[\s^*(t)].
\]
This problem is solved in two steps: firstly, one writes the Dynamic Programming Equation corresponding to the maximization problem for $J(u)$ {\em given} $m(t)$; then, one imposes that $m(t)$ is consistent with the master equation for the optimal process $\s^*(t)$. Denoting by $V(\s,t)$ the value function of the control problem of maximizing $J(u)$, the Dynamic Programming Equation reads, defining $\nabla V(\sigma,t) := V(-\s,t) - V(\s,t)$,
\begin{equation} \label{dp}
-\lambda V(\sigma,t)+\frac{\mu(\s,m(t))}{2}\left[\left[\nabla V(\sigma,t)\right]^+\right]^2+\frac{\partial V}{\partial t}(\sigma,t)+\sigma m(t)=0,
\end{equation}
and yields the optimal (feedback) control
\[
u^*(t) = \mu(\s,m(t)) \left[\nabla V(\sigma,t)\right]^+.
\]
A derivation of the Hamilton-Jacobi-Bellman equation (HJB), together with a classical verification argument, are postponed to Appendix A.
Substituting $u^*$ in \eqref{dyn}, one derives a differential equation for $m(t)$. It is convenient to write $\mu(\s,m)$ in the form $\mu(\s,m) = \s a(m) + b(m)$, and set $z(t) := \nabla V(1,t)$. By \eqref{dp} and \eqref{dyn} 
we obtain the following system of coupled equations:
\begin{equation}\label{3}
\left\{
\begin{array}{ll}
\dot{z}(t)&=\frac{b(m(t))}{2}z(t) |z(t)| +  \frac{a(m(t))}{2} z^2(t) +\lambda z(t)+2m(t)\\
\dot{m}(t)&=- (m(t) b(m(t)) + a(m(t))) |z(t)| \\
 &\qquad \qquad \  \qquad \qquad - (m(t) a(m(t)) + b(m(t))) z(t)  
\end{array}
\right . 
\end{equation}
% \begin{equation}\label{3}
% \left\{
% \begin{array}{lll}
% \dot{z}(t)=\frac{b(m(t))}{2}z(t) |z(t)| +  \frac{a(m(t))}{2} z^2(t) +\lambda z(t)+2m(t)\\
% \dot{m}(t)=- (m(t) b(m(t)) + a(m(t))) |z(t)| - (m(t) a(m(t)) + b(m(t))) z(t)  
% \end{array}
% \right . 
% \end{equation}
%We remark that the rigorous derivation of \eqref{3} as limit of the microscopic dynamic games does not follow from standard results, and is not the subject of the present work. Rigorous convergence results have been obtained recently for diffusion models  \cite{Go14,Ca15}. %More details in the heuristic derivation of \eqref{3} will be given in the Appendix.

Some remarks are needed concerning  equation \eqref{3}. It is relevant to note that  equation \eqref{3} should not be meant as an initial-value problem: only the initial $m(0)$, i.e. the initial information on agents' proportion, is assigned. On the other hand the value function in this problem is necessarily bounded, so only bounded solutions of \eqref{3} matter. Conversely, every bounded solution of \eqref{3} determines an equilibrium $u^*$ for the control problem associated to the functional given in \eqref{limmax}.

\section*{Baseline cases: constant mobility \& crowding effects}

In this section we consider two significant specifications of the model, for which we determine and characterize the bounded solutions of \eqref{3}. The proofs of the facts outlined below are postponed to Appendix B.

\subsection*{The constant mobility model} When $\mu(\s,m) = \mu = const$, \label{sec:standard} %In this model the mobility is constant.
equation \eqref{3} takes the form:
\begin{equation}\label{standard}
\left\{
\begin{array}{lll}
\dot{z}(t)=\frac{\mu}{2}z(t) |z(t)|  +\lambda z(t)+2m(t)\\
\dot{m}(t)=- \mu m(t) |z(t)| - \mu  z(t)  
\end{array}
\right . 
\end{equation}

We are interested in finding bounded solutions to \eqref{standard}. Note that $(z^*,m^*) = (0,0)$ is always an equilibrium. Moreover, two different regimes are detected under which the behavior of the system is completely different:

%\begin{theorem}\label{th:standard}
\begin{itemize}
\item[(a)] {\em Low mobility regime}: $\mu \leq \frac{\l^2}{8}$. For every $m(0) \in [-1,1]$ equation \eqref{standard} admits a unique bounded solution. For $m(0) \neq 0$ {\em consensus} occurs: $\lim_{t \ra +\infty} m(t) = {\mbox{sign}}(m(0)) \in \{-1,1\}$.
\item[(b)] {\em High mobility regime}: $\mu >  \frac{\l^2}{8}$. For $|m(0)| \neq 0$ sufficiently small there is more than one bounded solution to \eqref{standard}. All such solutions reach consensus ($\lim_{t \ra +\infty} m(t) \in \{\pm 1\}$), but exhibit a transient oscillatory regime, in which the orbits of the solutions spiral around $(0,0)$ before reaching consensus.
\end{itemize}
%\end{theorem}

\begin{figure}[h]
\begin{center}
\includegraphics[width=8.3cm,height=6.2cm]{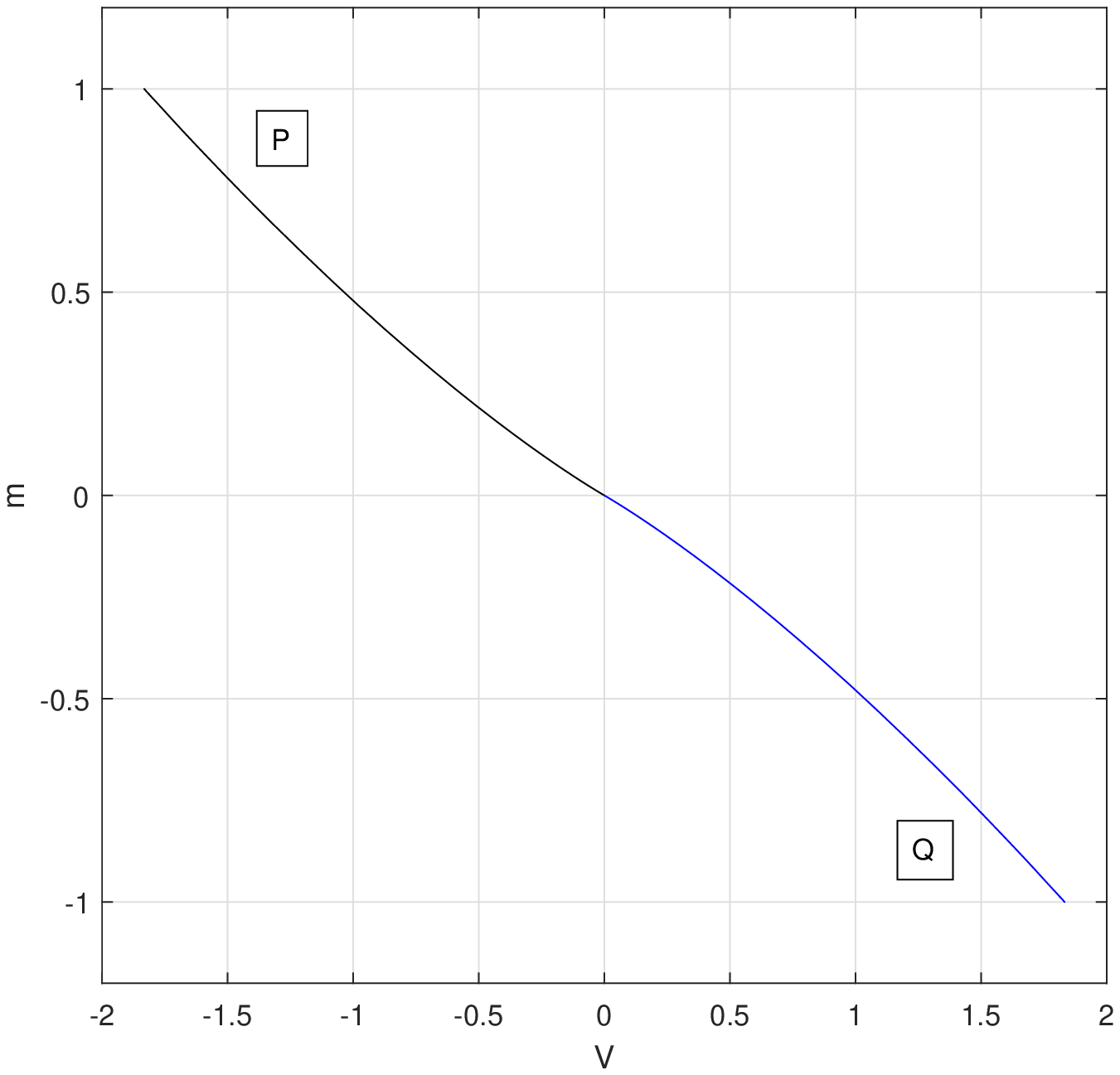}\\%{fig5}
\includegraphics[width=8.3cm,height=6.5cm]{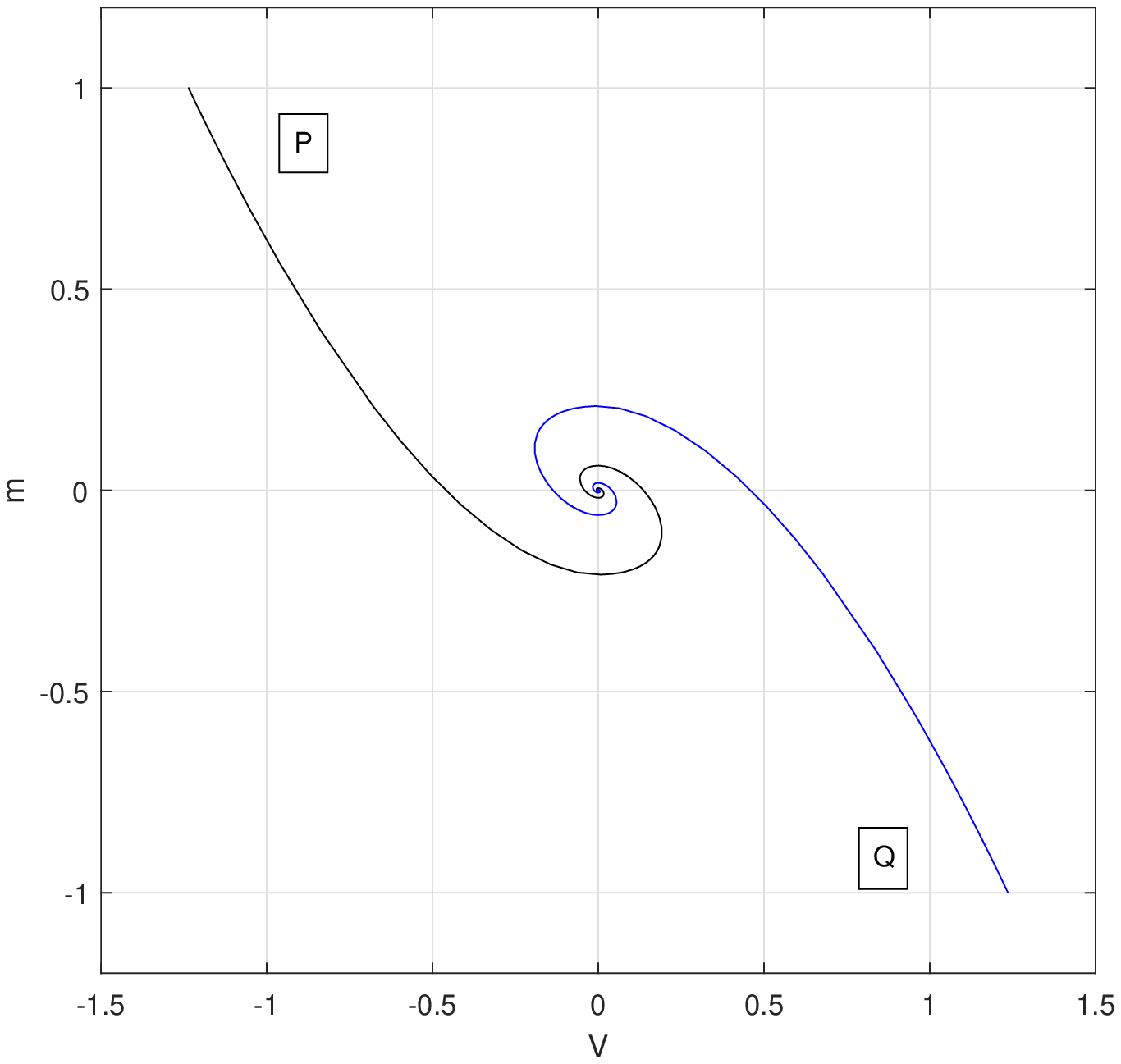}
\caption{Stable manifolds for the low mobility regime (top panel) and high mobility regime (bottom panel) in the case of constant mobility.} \label{fig3}
\end{center}
\end{figure}

In Fig. \ref{fig3} (top panel) we plot the stable manifolds related to the two fixed points $P$ and $Q$ of \eqref{standard}, different from the origin, for $\lambda=1$ and $\mu=0.1$. These values of the parameters fall under the low mobility regime. In Fig. \ref{fig3} (bottom panel), the values of the parameters are $\lambda=1$ and $\mu=1$. In this latter case, being under the high mobility regime, the manifolds are spiraling around the origin before reaching the consensus. Therefore, 
under this regime, the equilibrium control may be not unique: there are possibly multiple equilibrium controls leading to transient oscillating behavior.

\subsection*{Introducing crowding effects}\label{sec:crowding}

Here we set $\mu(\s,m) := \mu(1+ \e \s m)$, for some $\mu>0$ and $\e \in [0,1]$. Mobility is now asymmetric: changing state is more costly for an agent belonging to the minority: the cost to reinforce the majority increases the more the society is polarized. Put differently, the marginal cost to attract more people on the bandwagon is higher when the majority is more pronounced and, on the opposite, it is easier to loose some of them. Equation \eqref{3} becomes:
\begin{equation}\label{crowd}
\left\{
\begin{array}{lll}
\dot{z}(t)=\frac{\mu}{2}z(t) |z(t)|  +  \frac{\mu \e m}{2} z^2(t) +\lambda z(t)+2m(t)\\
\dot{m}(t)=- (1+\e)\mu m(t) |z(t)| - \mu(1+ \e m^2(t)) z(t) 
\end{array}
\right . 
\end{equation}

\begin{figure}[h]
\begin{center}
\includegraphics[width=7.5cm,height=6.5cm]{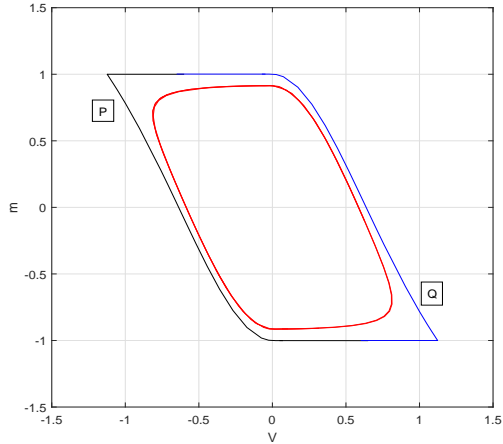}
\caption{Stable manifolds for $\mu=\hat \mu$ (blue and black curves) and a periodic orbit under the high mobility regime (red curve) in the case of crowding effects.} \label{fig6}
\end{center}
\end{figure}

Differently from the previous case, for certain values of the parameters,   an equilibrium control leading to {\em permanent} oscillatory behavior is detected. Indeed, a new threshold level $\hat{\mu}$, for the mobility parameter, with $\frac{\l^2}{8} < \hat{\mu} < +\infty$ appears. There are, therefore, three possible regimes:
%\begin{theorem} \label{th:crowd}
\begin{itemize}
\item[(a)] {\em Low mobility regime}: $\mu \leq \frac{\l^2}{8}$. For every $m(0) \in [-1,1]$ equation \eqref{crowd} admits a unique bounded solution. For $m(0) \neq 0$  {\em consensus} occurs: $\lim_{t \ra +\infty} m(t) = {\mbox{sign}}(m(0)) \in \{-1,1\}$.
\item[(b)] {\em  Moderate mobility regime}: $\frac{\l^2}{8}< \mu \leq \hat{\mu}$. For $|m(0)| \neq 0$ sufficiently small there is more than one bounded solution to \eqref{crowd}. All such solutions reach consensus ($\lim_{t \ra +\infty} m(t) \in \{\pm 1\}$), but, for $|m(0)|$ small enough, they exhibit a transient oscillatory regime, in which the orbits of the solutions spiral around $(0,0)$ before reaching consensus.
\item[(c)] {\em High mobility regime}: $\mu > \hat{\mu}$. For every $m(0) \in [-1,1]$ equation \eqref{crowd} admits two bounded solutions leading to consensus: $\lim_{t \ra +\infty} m(t)  \in \{-1,1\}$. Moreover \eqref{crowd} admits a unique non-constant periodic orbit: thus, for $|m(0)|$ sufficiently small, there are two periodic solutions which differ for a time shift.
\end{itemize}
In Fig. \ref{fig6} we plot the stable manifolds associated to the fixed points $P$ and $Q$ of \eqref{crowd}, for $\lambda=0.5$, $\epsilon=0.5$ and $\mu=\hat \mu=4.558$. Being $\mu$ exactly at its critical level, the two manifolds join at $P$ and $Q$. On the same graph, we have also depicted the periodic orbit obtained for $\lambda=0.5$, $\epsilon=0.5$ and with $\mu=4.6>\hat \mu$ (in red in the figure).

%\end{theorem}

What is the intuition behind the behavior of the population choosing an equilibrium corresponding to a periodic orbit as in the \emph{high mobility regime}? Without loss of genreality, suppose $m(t)$ is close to $1$ and  $\dot m(t)>0$ so that the society is on the way towards consensus. Because of the form of $\mu(\sigma, m)$, it is now very costly for a new agent to join the majority, whereas it is relatively less costly to abandon it. Moreover, note that the values of the parameters are such that the discount factor $\lambda$ (measuring impatience)  is small compared to $\mu$, so that the actual gain of being with the majority is relatively less important compared to future utility (where we can forecast a possible switch of the majority itself). Therefore, at some point, we reach a saturation value where \emph{it is not worth for a new agent to reach the majority}: the trade-off between the reward and the cost becomes negative. In some sense, agents have the perception that the run to climb the bandwagon is going to stop soon so that consensus will not be reached and therefore, sooner or later, the majority will change in favor of the party which is now the minority. Because of the cost structure proposed, it is now worth for agents in the majority to be the first to leave it (the cost is in fact lower if $m$ is high). In some sense, to be the first mover is more profitable. As a consequence $\dot m(t)$ changes its sign and we see a collapse to a negative value of $m$. Of course, when $m$ is close to $-1$, the simmetric situation happens and the society continues to oscillate forever.

\section*{Conclusions}
We have studied the dynamics of collective behavior in a society formed by a large population of conformist agents in charge to control for the effort to put in place in order to \emph{climb on the bandwagon}. Agents are remunerated by sticking with the majority but suffer a quadratic cost related to the control they put in place. Indeed, the binary state of the agents is controlled by progressively measurable controls determining the probability of switching the state. We determine the value function as the bounded solution to the Hamilton-Jacobi-Bellman equation and the relative optimal feedback control, which, it turns, also represents the Nash equilibrium of the associated mean-field game. 

We, then, solve in details two specifications of the model where the mobility parameter is either constant or crowed-dependent (in the sense that the cost to reach the majority increases the more polarized is the society). Interestingly enough, we find regimes of the parameter under which consensus is asymptotically reached and regimes under which bandwagoning is persistent  at the macro level: the systems is trapped in periodic solutions hence, it oscillates perpetually.

This paper sheds some light on a couple of open issues raised in the context of social interaction and collective behavior. Firstly, the possibility of detecting collective periodicity in complex social systems seemed to be unfeasible (see \cite{pnas16}). Moreover, differently from \cite{pnas16}, our agents solve a lifetime  optimization problem and can be, somewhat surprisingly, satisfied even with a permanently oscillating society. Secondly, our model provides an example of collective behavior on the short-time horizon which still accounts for diversity on the long-run: consensus is not necessarily the unique outcome of the society. With this respect, out model is in line with previous literature discussing long-term cultural diversity and short-term collective behavior (see \cite{va12}).

\subsection*{Appendix A. Derivation of the mean-field HJB equation}

\smallskip

Define
\[
V(\s,t) := \sup_{u}\E_{\s,t}\left[\int_t^{+\infty} e^{-\l(r-t)} R_u(r) dr \right],
\]
where
\[
R_u(t) :=\s(t) m(t) - \frac{1}{2\mu(\s(t), m(t))} u^2(t),
\]
$u(t)$ is the rate at which $\s(t)$ flips to $-\s(t)$, and $u = (u(t))_{t \geq 0}$ ranges over right continuous nonnegative closed loop controls, i.e. $u(t) = \varphi_u(t,\s(t))$ with $\varphi_u : [0,+\infty) \times \{-1,1\} \ra [0,+\infty]$ right continuous in $t$. Let
\[
J_{\s,t}(u) := \E_{\s,t}\left[\int_t^{+\infty} e^{-\l(r-t)} R_u(r) dr \right].
\]
If $u^*$ is an optimal control for \eqref{limmax}, then by the Bellman principle $V(\s,t) = J_{\s,t}(u^*)$ for all $\s,t$. For $t$ fixed and $h >0$, denote by $u^{h,\a}$ the control defined on $[t,+\infty)$ defined by
\[
u^{h,\a}(s) = \left\{ \begin{array}{ll} \a & \mbox{for } t \leq s < t+h \\ u^*(s) & \mbox{for } s \geq t+h, \end{array} \right.
\]
Observe that
%\begin{equation} \label{HJB1}
%\begin{split}
\begin{multline}\label{HJB1}
J_{\s,t}(u^{h,\a})   \\
=\E_{\s,t}\left[ \int_t^{t+h} R_{u^{h,\a}}(r) dr + e^{-\l h} V(\s(t+h),t+h) \right]  \\
=  h\left[\s(t) m(t) - \frac{1}{2\mu(\s(t), m(t))} \a^2 \right]   \\
  + \E_{\s,t}\left[e^{-\l h} V(\s(t+h),t+h) \right]+ o(h) .
% \end{split}
% \end{equation}
\end{multline}
Moreover
\begin{equation} \label{HJB2}
J_{\s,t}(u^{h,\a}) \geq V(\s,t)
\end{equation}
for every $\a$, while
\begin{equation} \label{HJB3}
J_{\s,t}(u^{h,u^*(t, \s)} )= V(\s,t) + o(h),
\end{equation}
where right continuity is used in this last estimate. It follows that
\begin{multline} \label{HJB4}
%\begin{split}
 \lim_{h \downarrow 0} \frac{\E_{\s,t}\left[e^{-\l h} V(\s(t+h),t+h)  - V(\s,t) \right]}{h}  \\
\geq \s(t) m(t) - \frac{1}{2\mu(\s(t), m(t))} \a^2
%\end{split}
\end{multline}
for every $\a \geq 0$, with the equality being attained at $\a = u^*(t, \s)$. By standard results on continuous time Markov chains
\begin{multline*}
\lim_{h \downarrow 0} \frac{\E_{\s,t}\left[e^{-\l h} V(\s(t+h),t+h) - V(\s,t) \right]}{h}  \\
= \frac{\partial V}{\partial t}(\sigma,t) +  \a \left[V(-\s,t) - V(\s,t)\right] - \l V(\s,t),
\end{multline*}
and \eqref{dp} follows.

We now show that if $(V(\s,t), m(t))$ solve \eqref{dp} coupled to the second equation in \eqref{3} and $V(\s,t)$ is {\em bounded}, then
\[
u^*(t) = \mu(\s,m(t)) \left[\nabla V(\sigma,t)\right]^+
\]
maximizes \eqref{limmax}. Note that the equation for $m(t)$ guarantees that $\E(\s^*(t)) = m(t)$, so $u^*$ is an equilibrium control. To show that $u^*$ maximizes \eqref{limmax} observe that
\begin{equation} \label{DP1}
\begin{split}
0 & = -\lambda V(\sigma,t)+u^*(t) \nabla V(\sigma,t) -\frac{1}{2 \mu(\s,m(t))}(u^*(t))^2 \\ & ~~~~+\frac{\partial V}{\partial t}(\sigma,t)+\sigma m(t)  \\
& = -\lambda V(\sigma,t) + \sup_a \left[ a \nabla V(\sigma,t) - \frac{1}{2 \mu(\s,m(t))}a^2 \right]\\ & ~~~~+ \frac{\partial V}{\partial t}(\sigma,t)+\sigma m(t) .
\end{split}
\end{equation}
Consider now an arbitrary feedback control $u$, and denote by $\s(t)$ the process with control $u$. A standard application of Ito's rule for Markov chains yields, for every $t>0$,
\begin{multline} \label{DP2}
\E\Bigg\{e^{-\l t} V(\s(t),t) - V(\s(0),0)  - \int_0^t \Big[ - \l e^{- \l s} V(\s(s),s) \\  + e^{- \l s} \frac{\partial V}{\partial s}(\s(s),s) - e^{-\l s} u(s) \nabla V(\s(s),s) \Big]ds \Bigg\} = 0
\end{multline}
Using \eqref{DP1}:
\begin{multline*}
u(s) \nabla V(\s(s),s) \leq \l V(\s(s),s) \\ +  \frac{1}{2 \mu(\s,m(s))} u^2(s) - \frac{\partial V}{\partial s}(\s(s),s) - \s(s) m(s),
\end{multline*}
which, inserted in \eqref{DP2} gives
\begin{multline} \label{DP3}
\E\Bigg\{e^{-\l t} V(\s(t),t) - V(\s(0),0) \\ + \int_0^t e^{- \l s} \left[\s(s) m(s) - \frac{1}{2 \mu(\s,m(s))} u^2(s) \right]ds \Bigg\} \leq 0,
\end{multline}
where equality is attained for $u = u^*$. Letting $t \rightarrow +\infty$ and using the boundedness of $V$, we obtain
\[
J(u) \leq \E[V(\s(0),0)] = J(u^*),
\]
and the proof is complete.

\subsection*{Appendix B. Derivation of other facts}\label{App}
%{\em Proof of Facts outlined in section \ref{sec:standard}}. 
\emph{Proof of Facts related to the constant mobility model}. 
We first observe that \eqref{standard}, besides the origin $O$, admits two other equilibria $P$ and $Q$, symmetric with respect to the origin: $\pm \left(({\sqrt{\l^2 + 4 \mu} - \l})/{\mu}, -1 \right)$. Linear analysis shows that $P$ and $Q$ are saddle points for all values of the parameters; the origin $O$ is linearly unstable:
\begin{itemize}
\item
for $\mu \leq \frac{\l^2}{8}$ it is repellent, i.e. the eigenvalues of the linearized system are both negative reals;
\item
for $\mu >  \frac{\l^2}{8}$ is an unstable spiral, i.e. the eigenvalues of the linearized system have both negative real part but nonzero imaginary part.
\end{itemize}
In order to perform a global analysis, we first consider the nullcline $\mathcal{N}$ given by the equation $\frac{\mu}{2}z|z| + \lambda z + 2m=0$. Off the nullcline, solutions to \eqref{standard} have trajectories that are locally graphs of a function $m = m(z)$. By implicit differentiation, assuming $(z,m) \in [0,+\infty) \times [-1,1]$, it turns out that $m''(z) > 0$ if and only if $\phi^-(z)<m<\phi^+(z)$, with
\[
\phi^{\pm}(z) = - \frac{z}{4} \left[  \l \mp \sqrt{\l^2 - 8 \mu + 6 \l \mu z + 4 \mu^2 z^2} \right].
\]
%\[
%\varphi^-(z) := - \frac{z}{4} \left[  \l + \sqrt{\l^2 - 8 \mu + 6 \l \mu z + 4 \mu^2 z^2} \right] < m < - \frac{z}{4} \left[  \l - \sqrt{\l^2 - 8 \mu + 6 \l \mu z + 4 \mu^2 z^2} \right]  =: \varphi^+(z).
%\]
For $(z,m) \in (-\infty,0) \times [-1,1]$, similar  convexity conditions are obtained by reflection w.r.t. the origin.  Consider the fixed point $Q$ and its stable manifold $\mathcal{M}_s$, i.e. the trajectory of a solution of \eqref{standard} converging to $Q$. 

\begin{figure}[h]
\begin{center}
\includegraphics[width=8cm,height=6.2cm]{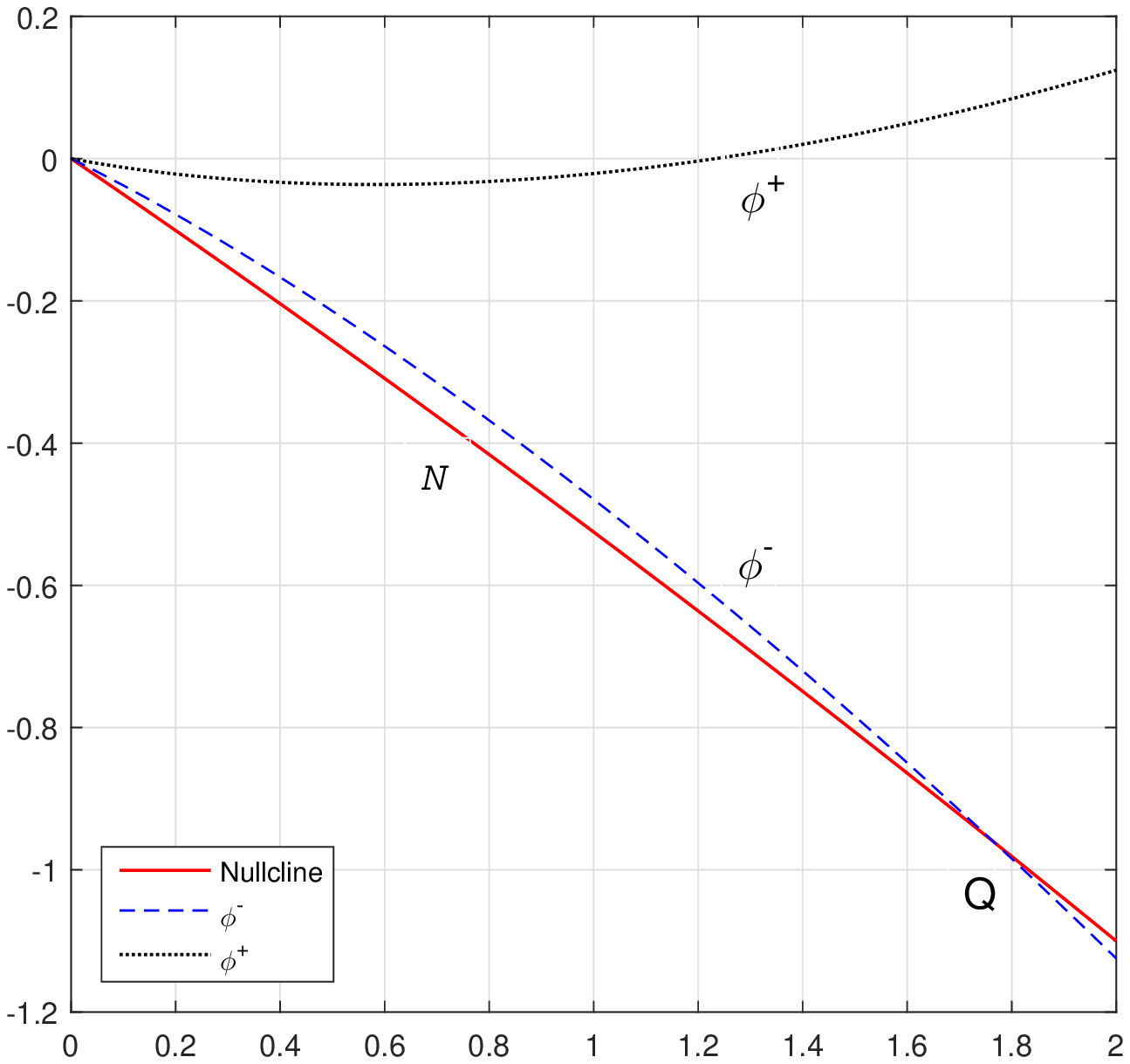}
\includegraphics[width=8cm,height=6.2cm]{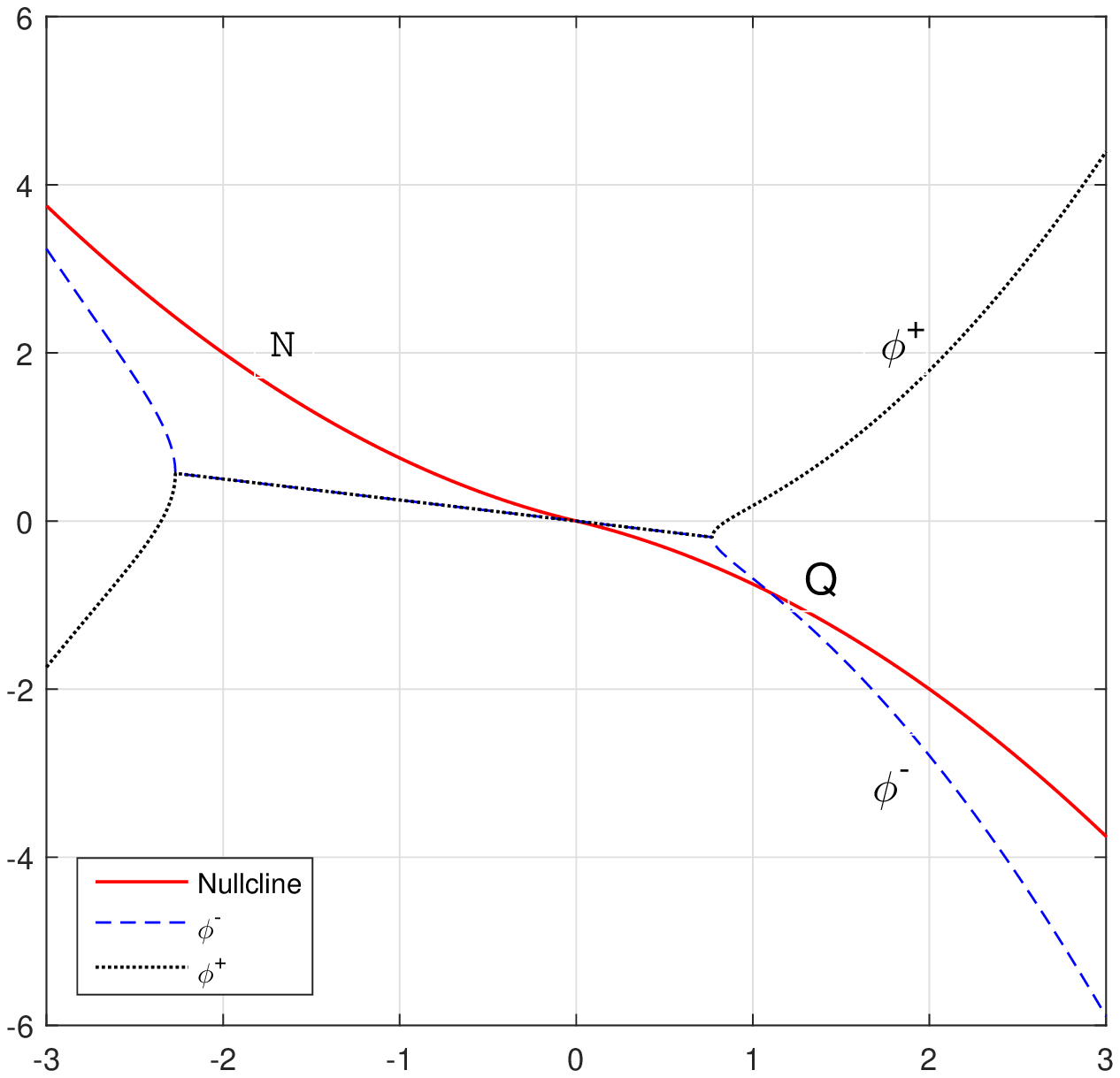}
\caption{Low mobility regime (top panel) and high mobility regime (bottom panel).} \label{fig1}
\end{center}
\end{figure}

\noindent
{\em Low mobility regime}: $\mu \leq \frac{\l^2}{8}$. In this case the graphs of $\phi^+$ and $\phi^-$ meet at the origin (see Fig. \ref{fig1}, top panel). Moreover, the graph of $\phi^-$ meets the nullcline $\mathcal{N}$ at the equilibrium point $Q$. A linear analysis at $Q$ and the study of the direction of the vector field of \eqref{standard} at the points of the graph of $\phi^-$ show that $\mathcal{M}_s$ is at the left of the graph of $\phi^-$. In particular $\mathcal{M}_s$ is concave, so it cannot intersect the nullcline $\mathcal{N}$, that can be intersected only vertically by a solution of \eqref{standard}. It follows that $\mathcal{M}_s$ is within the area between $\mathcal{N}$ and the graph of $\phi^-$. Since the origin is stable for the time-reversal of \eqref{standard}, necessarily $\mathcal{M}_s$ joins the origin with $Q$. Moreover, in the area between $\mathcal{N}$ and the graph of $\phi^-$, it easily checked that $\frac{dm}{dz} = \frac{\dot{m}}{\dot{z}} < 0$, so it is the graph of a strictly decreasing function. Thus, for every $m_0 \in (-1,0)$, there is a unique point of $\mathcal{M}_s$ with $m=m_0$, which is the starting point of a solution of \eqref{standard} converging to $Q$; in particular $m(t) \ra -1$ as $t \ra +\infty$. It is actually the only bounded solution starting from a point of the form $(m_0,z)$. This can be seen as follows. The point $(m_0,z)$, with $m_0<0$, cannot belong to the stable manifold of $P$, which is its image of $\mathcal{M}_s$ under reflection w.r.t the origin. Thus the solution starting from $(m_0,z)$ cannot converge to any fixed point. Moreover, since the divergence of the vector field driving \eqref{standard} is constantly equal to $\l >0$, then periodic orbits are not allowed. Thus, by the Poincar\'e-Bendixon Theorem, the solution starting from $(m_0,z)$ must be unbounded.

\noindent
{\em High mobility regime}: $\mu > \frac{\l^2}{8}$. In this case the graphs of $\phi^+$ and $\phi^-$ do not reach the origin (see Fig. \ref{fig1}, bottom panel). As in the low mobility regime, the stable manifold $\mathcal{M}_s$, as departing from $Q$, forms a concave curve between $\mathcal{N}$ and the graph of $\phi^-$. If we show that $\mathcal{M}_s$ gets arbitrarily close to the origin then the previous linear analysis implies that it must spiral around the origin, in particular it is not that graph of an injective function. 

Thus we are left to show that $\mathcal{M}_s$ gets arbitrarily close to the origin. This amounts to show that the solution $(\hat{z}(t),\hat{m}(t))$ of the time-reversed system starting from a point in $\mathcal{M}_s$ close to $Q$, converges to the origin as $t \ra +\infty$. Due to the spiraling around the origin, $(\hat{z}(t),\hat{m}(t))$ cannot converge to the origin following the graph of a monotone function. Thus it must intersect first the positive  $z$-axis and then the positive $m$ axis at some $m^*>0$. Suppose $m^*<1$. Note that $\mathcal{M}_s$ intersects the $m$-axis horizontally, so, again by convexity, after having touched $(0,m^*)$ it continues downward. Since $\mathcal{M}_s$, in the half-plane $z<0$ cannot touch the stable manifold of $P$, it follows it is trapped in a bounded region. Due to the absence of periodic orbits, necessarily $(\hat{z}(t),\hat{m}(t)) \ra (0,0)$ as $t \ra +\infty$.

Finally, we need to show that $m^*<1$.  By continuity from the low mobility regime, this is certainly true for $\mu - \frac{\l^2}{8}$ sufficiently small. If our claim is false, then there must be a value of $\mu$ for which $m^* = m^*(\mu) = 1$. In this situation, $\mathcal{M}_s$ continuous horizontally up to $P$. It follows that the union of $\mathcal{M}_s$ with the stable manifold of $P$ form a closed curve, tangent to the vector field driving \eqref{standard}; this is impossible by the Divergence Theorem.

\noindent
%{\em Sketch of the proof of Facts outlined in section \ref{sec:crowding}}. 
{\em Sketch of the proof of Facts related to the crowding effects model}. 
%For the following result we do not have a full proof yet. The numerical evidence appears solid, however, and leads to a rather clear picture of the bifurcations.

We first observe that  \eqref{crowd} has three equilibria: the origin $O$, whose linear properties are identical to those of the constant mobility model treated in the previous section, and the points $P$ and $Q$ with coordinates $\pm\left(-\frac{2}{\l},1 \right)$. Both $P$ and $Q$ are easily seen to be saddle points, for all values of the parameters. 
Similarly to  the constant mobility case,  the manifolds of $P$ and $Q$ can be proved to be  %: exactly as for the two phases of the standard model, these manifolds 
%are graphs of a 
monotone functions in the low mobility regime, while they spiral around the origin in the moderate mobility regime. What fails here is that the divergence of the driving vector field is not of constant sign, so that limit cycles cannot be ruled out. Although we do not have a full proof about the existence of a limit cycle, we provide clear evidence based on arguments derived by  numerical inspection.  Our analysis suggests that the $m$ coordinate of the  first intersection of the stable manifold of $Q$ with the $m$-axis is increasing in $\mu$, and it equals $1$ at some $\mu = \hat{\mu}$. Then, the manifold continues horizontally to reach $P$ (as depicted in Fig. \ref{fig6}). Thus, by symmetry, the two stable manifolds join to form a {\em separatrix}. By increasing $\mu$ further, a periodic orbit bifurcates from the separatrix through a homoclinic bifurcation.

\section*{Acknowledgments}{The authors acknowledge the financial support of Ca' Foscari University of Venice under the grant \emph{``Interactions in complex economic systems: innovation, contagion and crises''}. We are also grateful to Pietro Dindo, Tony He, Marco LiCalzi,  Paolo Pellizzari and other participants to seminars at Ca' Foscari, UTS Business School and attendants to the conference  MTNS 2016  for their precious comments. We are entirely responsible for all remaining errors. 
}

\bibliographystyle{acm} 
\bibliography{biblioDST17}

\begin{thebibliography}{10}

\bibitem{ak97social}
{\sc Akerlof, G.~A.}
\newblock Social distance and social decisions.
\newblock {\em Econometrica: Journal of the Econometric Society\/} (1997),
  1005--1027.

\bibitem{bafi17}
{\sc Bardi, M., and Fischer, M.}
\newblock On non-uniqueness and uniqueness of solutions in finite-horizon mean
  field games.
\newblock {\em arXiv preprint arXiv:1707.00628v2\/}, (2017).

\bibitem{JOTA16}
{\sc Bauso, D., Pesenti, R., and Tolotti, M.}
\newblock Opinion dynamics and stubbornness via multi-population mean-field
  games.
\newblock {\em Journal of Optimization Theory and Applications 170}, 1 (2016),
  266--293.

\bibitem{BeD}
{\sc Brock, W.~A., and Durlauf, S.~N.}
\newblock Discrete choice with social interactions.
\newblock {\em The Review of Economic Studies 68}, 2 (2001), 235--260.

\bibitem{Ca15}
{\sc Cardaliaguet, P., Delarue, F., Lasry, J.-M., and Lions, P.-L.}
\newblock The master equation and the convergence problem in mean field games.
\newblock {\em arXiv preprint arXiv:1509.02505\/} (2015).

\bibitem{cedafipe18}
{\sc Cecchin, A., Dai~Pra, P., Fischer, M., and Pelino, G.}
\newblock On the convergence problem in mean-field games: a two-state model
  without uniqueness.
\newblock In preparation (2018).

\bibitem{CeFi17}
{\sc Cecchin, A., and Fischer, M.}
\newblock Probabilistic approach to finite state mean field games.
\newblock {\em arXiv preprint arXiv:1704.00984\/} (2017).

\bibitem{CePe17}
{\sc Cecchin, A., and Pelino, G.}
\newblock Convergence, fluctuations and large deviations for finite state mean
  field games via the master equation.
\newblock {\em arXiv preprint arXiv:1707.01819\/} (2017).

\bibitem{ch97}
{\sc Challet, D., and Zhang, Y.-C.}
\newblock Emergence of cooperation and organization in an evolutionary game.
\newblock {\em Physica A: Statistical Mechanics and its Applications 246}, 3-4
  (1997), 407--418.

\bibitem{dp15}
{\sc Collet, F., Dai~Pra, P., and Formentin, M.}
\newblock Collective periodicity in mean-field models of cooperative behavior.
\newblock {\em NoDEA - Nonlinear Differential Equations and Applications
  22(5)\/} (2015), 1461--1482.

\bibitem{dp13}
{\sc Dai~Pra, P., Fischer, M., and Regoli, D.}
\newblock A curie-weiss model with dissipation.
\newblock {\em Journal of Statistical Physics 152(1)\/} (2013), 37--53.

\bibitem{dit17}
{\sc Ditlevsen, S., and L{\"o}cherbach, E.}
\newblock Multi-class oscillating systems of interacting neurons.
\newblock {\em Stochastic Processes and their Applications 127}, 6 (2017),
  1840--1869.

\bibitem{ep02}
{\sc Epstein, J.~M.}
\newblock Modeling civil violence: An agent-based computational approach.
\newblock {\em Proceedings of the National Academy of Sciences 99}, suppl 3
  (2002), 7243--7250.

\bibitem{fo74}
{\sc F{\"o}llmer, H.}
\newblock Random economies with many interacting agents.
\newblock {\em Journal of mathematical economics 1}, 1 (1974), 51--62.

\bibitem{basar2016}
{\sc Gao, Z., Chen, X., Liu, J., and Basar, T.}
\newblock Periodic behavior of a diffusion model over directed graphs.
\newblock {\em Decision and Control (CDC), 2016 IEEE 55th Conference on.
  IEEE\/} (2016).

\bibitem{Go14}
{\sc Gomes, D.~A., et~al.}
\newblock Mean field games models: a brief survey.
\newblock {\em Dynamic Games and Applications 4}, 2 (2014), 110--154.

\bibitem{Gr}
{\sc Granovetter, M.}
\newblock Threshold models of collective behavior.
\newblock {\em American journal of sociology 83}, 6 (1978), 1420--1443.

\bibitem{Gu11}
{\sc Gu\'{e}ant, O., Lasry, J.-M., and Lions, P.-L.}
\newblock Mean field games and applications.
\newblock In {\em Paris-Princeton lectures on mathematical finance 2010},
  vol.~2003 of {\em Lecture Notes in Mathematics}. Springer Berlin Heidelberg,
  2011, pp.~205--266.

\bibitem{ho06}
{\sc Horst, U., and Scheinkman, J.~A.}
\newblock Equilibria in systems of social interactions.
\newblock {\em Journal of Economic Theory 130}, 1 (2006), 44--77.

\bibitem{hu03}
{\sc Huang, M., Caines, P.~E., and Malham{\'e}, R.~P.}
\newblock Individual and mass behaviour in large population stochastic wireless
  power control problems: centralized and nash equilibrium solutions.
\newblock In {\em Decision and Control, 2003. Proceedings. 42nd IEEE Conference
  on\/} (2003), vol.~1, IEEE, pp.~98--103.

\bibitem{caines2005}
{\sc Huang, M., Malham\'e, R.~P., and Caines, P.~E.}
\newblock Nash equilibria for large-population linear stochastic systems of
  weakly coupled agents.
\newblock {\em Analysis, Control and Optimization of Complex Dynamic Systems\/}
  (2005), 215--252.

\bibitem{Ca06}
{\sc Huang, M., Malham{\'e}, R.~P., Caines, P.~E., et~al.}
\newblock Large population stochastic dynamic games: closed-loop mckean-vlasov
  systems and the nash certainty equivalence principle.
\newblock {\em Communications in Information \& Systems 6}, 3 (2006), 221--252.

\bibitem{kalai04}
{\sc Kalai, E.}
\newblock Large robust games.
\newblock {\em Econometrica 72}, 6 (2004), 1631--1665.

\bibitem{Li06}
{\sc Lasry, J.-M., and Lions, P.-L.}
\newblock Jeux \`{a} champ moyen. i -- le cas stationnaire.
\newblock {\em Comptes Rendus Mathematique 343(9)\/} (2006), 619--625.

\bibitem{Li04}
{\sc Lindner, B., Garc{\i}a-Ojalvo, J., Neiman, A., and Schimansky-Geier, L.}
\newblock Effects of noise in excitable systems.
\newblock {\em Physics reports 392}, 6 (2004), 321--424.

\bibitem{mo10}
{\sc Montanari, A., and Saberi, A.}
\newblock The spread of innovations in social networks.
\newblock {\em Proceedings of the National Academy of Sciences 107}, 47 (2010),
  20196--20201.

\bibitem{pnas16}
{\sc Ramazi, P., Riehl, J., and Cao, M.}
\newblock Networks of conforming or nonconforming individuals tend to reach
  satisfactory decisions.
\newblock {\em Proceedings of the National Academy of Sciences 113}, 46 (2016),
  12985--12990.

\bibitem{sh78_micro}
{\sc Schelling, T.}
\newblock Micromotives and macrobehavior.
\newblock {\em Norton 372\/} (1978), 373.

\bibitem{sc71}
{\sc Schelling, T.~C.}
\newblock Dynamic models of segregation.
\newblock {\em Journal of mathematical sociology 1}, 2 (1971), 143--186.

\bibitem{sc73_helmets}
{\sc Schelling, T.~C.}
\newblock Hockey helmets, concealed weapons, and daylight saving: A study of
  binary choices with externalities.
\newblock {\em Journal of Conflict resolution 17}, 3 (1973), 381--428.

\bibitem{Ieee13}
{\sc Stella, L., Bagagiolo, F., Bauso, D., and Como, G.}
\newblock Opinion dynamics and stubbornness through mean-field games.
\newblock In {\em Decision and Control (CDC), 2013 IEEE 52nd Annual Conference
  on\/} (2013), IEEE, pp.~2519--2524.

\bibitem{tou14}
{\sc Touboul, J.}
\newblock The hipster effect: When anticonformists all look the same.
\newblock {\em arXiv preprint arXiv:1410.8001\/} (2014).

\bibitem{va12}
{\sc Valori, L., Picciolo, F., Allansdottir, A., and Garlaschelli, D.}
\newblock Reconciling long-term cultural diversity and short-term collective
  social behavior.
\newblock {\em Proceedings of the National Academy of Sciences 109}, 4 (2012),
  1068--1073.

\bibitem{yang2016}
{\sc Yang, T., Meng, Z., Dimarogonas, D.~V., and Johansson, K.~H.}
\newblock Periodic behaviors for discrete-time second-order multiagent systems
  with input saturation constraints.
\newblock {\em IEEE Transactions on Circuits and Systems II: Express Briefs
  63}, 7 (2016), 663--667.

\bibitem{yo11}
{\sc Young, H.~P.}
\newblock The dynamics of social innovation.
\newblock {\em Proceedings of the National Academy of Sciences 108}, Supplement
  4 (2011), 21285--21291.

\end{thebibliography}

\end{document}